\documentclass{amsart}
\usepackage{amssymb,latexsym,amsmath,amsfonts,hhline,verbatim,color}
\usepackage{pdfsync}
\input xy
\xyoption{all}

\newcommand{\cal}{\mathcal}

\makeatletter
\renewcommand{\subsection}{\@startsection{subsection}{2}{0mm}{-2mm}{-2mm}{\bf\normalsize}}

\def\sbsnt#1{\subsection{#1}}
\makeatother

\newtheorem{formula}{}[section]
\newtheorem{definition}[formula]{Definition}
\newtheorem{corollary}[formula]{Corollary}
\newtheorem{remark}[formula]{Remark}
\newtheorem{lemma}[formula]{Lemma}
\newtheorem{theorem}[formula]{Theorem}

\def\thrm{\begin{theorem}}
\def\thrml#1{\begin{theorem}\label{#1}}
\def\ethrm{\end{theorem}}
\def\rmrk{\begin{remark}}
\def\rmrkl#1{\begin{remark}\label{#1}}
\def\ermrk{\end{remark}}
\def\dfntn{\begin{definition}}
\def\dfntnl#1{\begin{definition}\label{#1}}
\def\edfntn{\end{definition}}
\def\nmrt{\begin{enumerate}}
\def\enmrt{\end{enumerate}}
\def\tm#1{\item[{\rm (#1)}]}
\def\qtn{\begin{equation}}
\def\qtnl#1{\begin{equation}\label{#1}}
\def\eqtn{\end{equation}}
\def\lmm{\begin{lemma}}
\def\lmml#1{\begin{lemma}\label{#1}}
\def\elmm{\end{lemma}}
\def\crllr{\begin{corollary}}
\def\crllrl#1{\begin{corollary}\label{#1}}
\def\ecrllr{\end{corollary}}
\def\css{\begin{cases}}
\def\ecss{\end{cases}}

\def\proof{\noindent{\bf Proof}.\ }

\def\cG{{\cal G}}
\def\cL{{\cal L}}
\def\cK{{\cal K}}

\def\cS{{\cal S}}

\def\cX{{\cal X}}
\def\cY{{\cal Y}}

\def\mF{{\mathbb F}}

\def\SG{{\scriptscriptstyle\Gamma}}

\DeclareMathOperator{\aut}{Aut}

\DeclareMathOperator{\diag}{Diag}

\DeclareMathOperator{\Fix}{Fix}
\DeclareMathOperator{\fix}{fix}

\DeclareMathOperator{\id}{id}

\DeclareMathOperator{\inv}{Inv}

\DeclareMathOperator{\iso}{Iso}

\DeclareMathOperator{\orb}{Orb}

\DeclareMathOperator{\PSL}{PSL}

\DeclareMathOperator{\SL}{SL}

\DeclareMathOperator{\syl}{Syl}
\DeclareMathOperator{\sym}{Sym}
\DeclareMathOperator{\WL}{WL}

\def\bull{\hfill\vrule height .9ex width .8ex depth -.1ex\medskip}

\def\qaq{\quad\text{and}\quad}
\def\qoq{\quad\text{or}\quad}
\def\ov{\overline}

\def\wt{\widetilde}

\begin{document}
\title{Cartan coherent configurations}
\author{Ilia Ponomarenko}
\address{St.Petersburg Department of the Steklov Mathematical
Institute, St.Petersburg, Russia}
\email{inp@pdmi.ras.ru}
\author{Andrey Vasil'ev}
\thanks{The work of the first and the second authors was partially
supported respectively, by the grant RFBR 
No.~14-01-00156 and RFFI Grant No.~13-01-00505}
\address{Sobolev Institute of Mathematics, Novosibirsk, Russia\vspace{-7pt}}
\address{Novosibirsk State University, Novosibirsk, Russia}
\email{vasand@math.nsc.ru}

\begin{abstract}
The Cartan scheme $\cX$ of a finite group $G$ with a $(B,N)$-pair is defined
to be the coherent configuration associated with the action of $G$ on the
right cosets of the Cartan subgroup $B\cap N$ by the right multiplications.
It is proved that if $G$ is a simple group of Lie type, then asymptotically,
the coherent configuration~$\cX$ is $2$-separable, i.e., the array of 
$2$-dimensional intersection
numbers determines $\cX$ up to isomorphism. It is also proved that in this
case, the base number of~$\cX$ equals~$2$. This enables us
to construct a polynomial-time algorithm for recognizing the Cartan schemes
when the rank of $G$ and order of the underlying field are sufficiently
large. One of the key points in the proof of the main results is a new
sufficient condition for an arbitrary homogeneous coherent configuration to 
be $2$-separable.
\end{abstract}

\maketitle

\section{Introduction}

A well-known general problem in algebraic combinatorics is to characterize
an association scheme $\cX$ up to isomorphism by a certain set of
parameters~\cite{BI}. A lot of such characterizations are known when $\cX$ 
is the association scheme of a classical distance regular graph~\cite{BCN}.
In most cases, the parameters can be chosen as a part of the intersection 
number array of~$\cX$. However in general, even the whole array does not
determine the scheme $\cX$ up to isomorphism. Therefore, it makes sense 
to consider the $m$-dimensional intersection numbers, $m\ge 1$,
introduced in~\cite{EP00} for an arbitrary coherent configuration 
(for $m=1$, these numbers are ordinary intersection numbers; the exact
definitions can be found in Section~\ref{110415}). It was proved in~\cite{EP00}
that every  Johnson, Hamming or Grassmann scheme is $2$-separable, i.e.,
is determined up to isomorphism by the array of $2$-dimensional 
intersection numbers.\medskip

In a recent paper \cite{APM}, a generalized notion of distance regularity 
in buildings was introduced. According to \cite{{Zi2}}, there is a natural 
1-1 correspondence between the class of all buildings and the class of 
special homogeneous coherent configurations called the Coxeter schemes
(see also~\cite[Chapter~12]{Zi1}).
In this language, the Tits theorem on spherical buildings says 
that if  $\cX$ is
a  finite Coxeter scheme of rank at least~$3$, then there exists
a group $G$ acting on a set $\Omega$ such that
\qtnl{290116a}
\cX=\inv(G,\Omega)
\eqtn
where $\inv(G,\Omega)$ is the coherent configuration of $G$, i.e., 
the pair $(\Omega,S)$ with $S=\orb(G,\Omega\times\Omega)$. Moreover,
in this case, $G$ is a group with a $(B,N)$-pair. 
Thus, a characterization of the coherent configuration~\eqref{290116a}
with such $G$ by the $m$-dim intersection numbers with small $m$ could 
be consider as a generalization of the above mentioned results on the
association  schemes of classical distance regular graphs to the
noncommutative case.\medskip

In the present paper, we are interested in coherent
configurations~\eqref{290116a} in the case when $G$ is a finite
group with a $(B,N)$-pair and $\cX$ is a {\it Cartan scheme of $G$} 
in the following sense.

\dfntn\label{060216t}
The Cartan scheme of $G$ with respect to $(B,N)$ is defined to be
the coherent configuration~\eqref{290116a}, where $\Omega=G/H$ consists
of the right cosets of the Cartan subgroup $H=B\cap N$ and $G$
acts on $\Omega$ by right multiplications.
\edfntn

Note that the permutation group induced by the action of $G$ is
transitive and the stabilizer of the point $\{H\}$ coincides with
$H$. In a Coxeter scheme of rank at least~$3$, the point
stabilizer equals~$B$. Therefore, every  Coxeter scheme is a quotient
of a suitable Cartan scheme.\medskip

The separability problem~\cite{EP09} consists in finding the
smallest $m$, for which a coherent configuration is $m$-separable.
The separability problem (in particular, for a Cartan scheme)
is easy to solve if the group $H$ is trivial. Indeed, in
this case, the permutation group induced by $G$ is regular and the
corresponding coherent configuration is 1-separable. The following theorem
gives a partial solution to the separability problem for Cartan 
schemes when $G$ is a finite simple group of Lie type, and hence with a (B,N)-pair.
In what follows, we denote by $\cL$ the class 
of all simple groups of Lie type including all exceptional groups
and all classical groups $\Phi(l,q)$, for which $l\ge l_0$ and $q\ge al$,
where the values of $l_0$ and $a$ are given in the last two columns 
of Table~\ref{t2} at page~\pageref{010216a}.

\thrml{030415a}
The Cartan scheme $\cX$ of every finite simple group~$G\in\cL$ 
is $2$-separable.
\ethrm

As a byproduct of the proof of Theorem~\ref{030415a}, we are able
to estimate the base number of a Cartan scheme satisying the
hypothesis of this theorem (as to the exact definition,
we refer to Subsection~\ref{300116a}; see also~\cite{EP09}
and~\cite[Sec.~5]{BC}, where the base number was called the
EP-dimension). The base number $b(\cX)$ of a coherent
configuration $\cX$ can be thought 
as a combinatorial analog of the base number of a permutation group,
which is the minimal number of points such that only identity of 
the group leave each of them fixed. In fact, for the coherent configuration~\eqref{290116a},
we have
\qtnl{200615a}
b(G)\le b(\cX),
\eqtn
where $b(G)$ is the base number of the permutation group induced
by~$G$. Moreover, in this case, obviously, $b(G)=1$ if and only 
if $b(\cX)=1$. In general, $b(G)$ can be much smaller than $b(\cX)$.
The following theorem shows that this does not happen for the
Cartan schemes in question.

\thrml{300116c}
Let $\cX$ be the Cartan scheme of a group $G\in\cL$. Then 
$b(\cX)\le 2$ and $b(\cX)=1$ if and only if the group $H$ is trivial. 
\ethrm

Let us deduce Theorems~\ref{030415a} and~\ref{300116c} from the
results, which proofs occupy the most part of the paper.
Let $\cX$ be the Cartan scheme of a group $G\in\cL$. Denote 
by $c$ and $k$  the indistinguishing number and the maximum valency 
of $\cX$, respectively.\footnote{In the complete colored graph
representing $\cX$,	$k$ is the maximum number of the monochrome 
arcs incident to a vertex, and $c$ is the  maximum number of 
triangles with fixed base, the other two sides of which are 
monochrome arcs.} Translating these invariants into the
group-theoretic language, we prove in Theorem~\ref{310116b} that
in our case
\qtnl{120415a}
2c(k-1)<n
\eqtn
where $n=|\Omega|$. The proof of this inequality  forms the 
group-theoretic part of the whole proof. The
combinatorial part of the proof  is
to analyze one point extension of a homogeneous coherent
configuration, for which inequality~\eqref{120415a} holds;
here the point extension can be thought as a combinatorial analog 
of the point stabilizer of a permutation group. In this way, we prove
Theorem~\ref{291114a}, which implies that the one point extension
of $\cX$ is $1$-regular (see Subsection~\ref{300116a}). 
It immediately follows
that $b(\cX)\le 2$, which proves Theorem~\ref{300116c}. Finally, 
in view of Theorem~\ref{291114a}, the $2$-separability of $\cX$
is a direct consequence of Theorem~\ref{310116c}
obtained by a combination of two results in~\cite{EP09}, so
Theorem~\ref{030415a} holds.\medskip

When the rank of a simple  group $G$ of Lie type is small,
inequality~\eqref{120415a} does  not generally hold, but the 
statements of Theorems~\ref{030415a} and~\ref{300116c} may still 
be true. For example, this happens in the following case.

\thrml{120714a}
Let $\cX$ be the Cartan scheme of the group $\PSL(2,q)$, where $q>3$.
Then $\cX$ is $2$-separable and $b(\cX)=2$.
\ethrm\medskip

We believe that the Cartan scheme of every simple group of Lie type
is $2$-separable. Moreover, as in the case of classical distance
regular graphs, it might be that in most cases 
such a scheme is $1$-separable, i.e., is determined up to isomorphism
by the intersection numbers. In this way, one could probably use
more subtle results on the structure of  finite simple groups
and a combinatorial technique in spirit of~\cite{MP12a}.\medskip

From the computational point of view, Theorems~\ref{030415a}
and~\ref{300116c} can be used for testing isomorphism and 
recognizing the Cartan schemes satisfying the hypothesis
of Theorem~\ref{030415a}. For this aim, it  is convenient to 
represent a coherent configuration $(\Omega,S)$ as a 
complete colored graph with vertex set $\Omega$, in which the 
color classes of arcs 
coincide with relations of the set $S$ (the vertex colors match 
the colors of the loops). It is assumed that isomorphisms
of such colored graphs preserve the colors.

\thrml{310116e}
Let $\cG_n$ (resp. $\cK_n$) be the class of all  colored
graphs (resp. the colored graphs of Cartan schemes of the
groups in $\cL$) with $n$ vertices.
Then the  following problems can be solved in polynomial time in $n$:
\nmrt 
\tm{1} given $D\in\cG_n$, test whether $D\in\cK_n$, and 
(if so) find the corresponding groups $G$, $B$, and $N$;
\tm{2} given $D\in\cK_n$ and $D'\in\cG_n$, find the set $\iso(D,D')$.
\enmrt
\ethrm\medskip

To make the paper possibly self-contained, we cite the basics of
coherent configurations in Section~\ref{110415}.
Theorems~\ref{291114a} and~\ref{310116b}, from which we have deduced
Theorems~\ref{030415a} and~\ref{300116c}, are proved in
Section~\ref{250415a} and Sections~\ref{310116a},\ref{070715a},
respectively. Finally, Theorems~\ref{120714a} and~\ref{310116e}
are proved in Sections~\ref{240714z} and~\ref{210615b}, respectively.
\medskip

{\bf Notation.}
Throughout the paper, $\Omega$ denotes a finite set.

The diagonal of the Cartesian product $\Omega\times\Omega$ is denoted by~$1_\Omega$; for any $\alpha\in\Omega$
we set $1_\alpha=1_{\{\alpha\}}$.

For a relation $r\subset\Omega\times\Omega$, we set $r^*=\{(\beta,\alpha):\ (\alpha,\beta)\in r\}$ and
$\alpha r=\{\beta\in\Omega:\ (\alpha,\beta)\in r\}$ for all $\alpha\in\Omega$.

For $S\in 2^{\Omega^2}$, we denote by $S^\cup$ the set of all unions of the elements of $S$, and put $S^*=\{s^*:\ s\in S\}$
and $\alpha S=\cup_{s\in S}\alpha s$, where $\alpha\in\Omega$.

For $g\in\sym(\Omega)$, we set $\Fix(g)=\{\alpha\in\Omega:\ \alpha^g=\alpha\}$; in particular, if  $\chi$ is the permutation
character of a group $G\le\sym(\Omega)$, then
$\chi(g)=|\Fix(g)|$ for all $g\in G$.

The identity of a group $G$ is denoted by $e$; the set of non-identity elements in $G$ is denoted by  $G^\#$.

\section{Coherent configurations}\label{110415}
\sbsnt{Main definitions.}\label{230415b}
Let $\Omega$ be a finite set, and let $S$ be a partition of~$\Omega\times\Omega$. The pair $\cX=(\Omega,S)$ is called
a {\it coherent configuration} on $\Omega$ if $1_\Omega\in S^\cup$, $S^*=S$, and given $r,s,t\in S$, the number
$$
c_{rs}^t=|\alpha r\cap\beta s^*|
$$
does not depend on the choice of $(\alpha,\beta)\in t$. The elements of $\Omega$, $S$, and $S^\cup$ are
called the {\it points}, {\it basis relations}, and {\it relations} of~$\cX$, respectively.
The numbers $|\Omega|$, $|S|$, and $c_{rs}^t$ are called the {\it degree}, {\it rank}, and {\it intersection numbers} of~$\cX$. 
The basis relation containing the pair $(\alpha,\beta)\in\Omega\times\Omega$ is denoted by $r(\alpha,\beta)$.\medskip

The point set $\Omega$ is a disjoint union of {\it fibers}, i.e., the sets $\Gamma\subseteq\Omega$, for which $1_\SG\in S$
For any basis relation $r\in S$, there exist uniquely determined fibers $\Gamma$ and $\Delta$
such that $r\subseteq\Gamma\times\Delta$. Moreover, the number $|\gamma r|=c_{rr^*}^t$ with $t=1_\SG$,
does not depend on the choice of $\gamma\in\Gamma$. This number is called the {\it valency} of~$r$ and denoted~$n_r$. The maximum of all
valences is denoted by $k=k(\cX)$.\medskip

A point $\alpha\in\Omega$ of the coherent conﬁguration $\cX$ is called {\it regular} if
$$
|\alpha r|\le 1\quad\text{for all}\ \,r\in S. 
$$
One can see that the set of all regular points is the union of fibers. If this set is not empty, then the coherent conﬁguration 
$\cX$ is said to be {\it 1-regular}.\medskip

The coherent configuration $\cX$ is said to be {\it homogeneous}
if $1_\Omega\in S$. In this case,
$n_r=n_{r^*}=|\alpha r|$ for all $r\in S$ and $\alpha\in\Omega$. Moreover, the relations
\qtnl{150410a}
c_{r^*s^*}^{t^*}=c_{sr}^t\quad\text{and}\quad
n_tc_{rs}^{t^*}=n_rc_{st}^{r^*}=n_sc_{tr}^{s^*}
\eqtn
hold for all $r,s,t\in S$. We observe that in the homogeneous case, a coherent configuration is 1-regular if and only if it is a thin
scheme in the sense of~\cite{Zi1}.

\sbsnt{Point extensions and  the base number.}\label{300116a}
There is a natural partial order\, $\le$\, on the set of all coherent configurations on the same set.
Namely, given two coherent configurations $\cX=(\Omega,S)$ and
$\cX'=(\Omega,S')$, we set
$$
\cX\le\cX'\ \Leftrightarrow\ S^\cup\subseteq (S')^\cup.
$$
 The minimal and maximal elements with respect to this ordering are the {\it trivial} and {\it complete} coherent
configurations: the basis relations of the former one are the reflexive relation $1_\Omega$ and
(if $n>1$) its complement in $\Omega\times\Omega$, whereas the basis relations of the latter one are singletons.\medskip

Given two coherent configurations $\cX_1$ and $\cX_2$ on $\Omega$, there is a uniquely determined coherent configuration
$\cX_1\cap\cX_2$ also on $\Omega$, the relation set of which is $(S_1)^\cup\cap(S_2)^\cup$, where $S_i$ is the set of basis
relations of~$\cX_i$, $i=1,2$. This enables us to define the {\it point extension} $\cX_{\alpha,\beta,\ldots}$ of a coherent configuration
$\cX=(\Omega,S)$ with respect to the points $\alpha,\beta,\ldots\,\in\Omega$ as follows:
$$
\cX_{\alpha,\beta,\ldots}=\bigcap_{\cY:\ S\subseteq T^\cup,1_\alpha,1_\beta,\ldots\in T^\cup}\cY,
$$
where $\cY$ is the coherent configuration $(\Omega,T)$. In other words,
$\cX_{\alpha,\beta,\ldots}$ can be defined as the smallest coherent configuration 
on $\Omega$ that is larger than or equal to $\cX$ and has singletons
$\{\alpha\},\{\beta\},\ldots$ as fibers. This configuration can also be considered as the refinement of the color
graph associated  with $\cX$, in which the points of $\alpha,\beta,\ldots$ are colored in distinguished new colors. In
particular, the extension can be efficiently constructed by the
Weisfeiler-Leman algorithm (see Section~\ref{210615b}).

\dfntn
A set $\{\alpha,\beta,\ldots\}\subseteq\Omega$   is a {\it base} of the coherent configuration $\cX$ if  the extension
$\cX_{\alpha,\beta,\ldots}$  with respect to the points $\alpha,\beta,\ldots$
is complete; the smallest cardinality of a base is called the {\it base number} of $\cX$ and denoted by $b(\cX)$.
\edfntn

It is easily seen that $0\le b(\cX)\le n-1$, where $n=|\Omega|$, and the equalities are attained for the complete and trivial coherent
configurations on $\Omega$, respectively. It is also obvious that $b(\cX)\le 1$, whenever the 
coherent configuration $\cX$ is 1-regular.

\sbsnt{Coherent configurations and permutation groups.}\label{030216f}
Two coherent configurations $\cX=(\Omega,S)$ and $\cX'=(\Omega',S')$
are called {\it isomorphic} if there exists a bijection
$f:\Omega\to\Omega'$ such that the relation
$s^f=\{(\alpha^f,\beta^f):\ (\alpha,\beta)\in s\}$ belongs to $S'$
for all $s\in S$. The bijection $f$ is called an {\it isomorphism} 
from $\cX$ onto $\cX'$; the set of all of them is denoted 
by $\iso(\cX,\cX')$. The group $\iso(\cX,\cX)$ contains a normal subgroup
$$
\aut(\cX)=\{f\in\sym(\Omega):\ s^f=s,\ s\in S\}
$$
called the {\it automorphism group} of~$\cX$.\medskip 

Let $G\le\sym(\Omega)$ be a permutation group, and let $S$ be the 
set of orbits of the coordinatewise action of $G$
on~$\Omega\times\Omega$. Then 
$$
\inv(G)=\inv(G,\Omega)=(\Omega,S)
$$
is a coherent configuration called the {\it coherent 
configuration of~$G$}. It is homogeneous if and only if the group 
$G$ is transitive. From \cite[Corollary~3.4]{EP09},
it follows that a coherent configuration $\cX$ is 1-regular if and only if $\cX=\inv(G)$, where $G$ is a permutation group having 
a faithful regular orbit.\medskip

Let $G\le\sym(\Omega)$ be a transitive group, $H=G_\alpha$ the 
stabilizer of a point $\alpha$ in~$G$,
and $\cX=\inv(G)$ the coherent configuration  of~$G$. Then given a basis relation $s\in S$, one can form
the set
\qtnl{140714u}
D_s=\{g\in G:\ (\alpha,\alpha^g)\in s\},
\eqtn
which is, in fact, a double $H$-coset. It is  well known that the
mapping $s\mapsto D_s$ is a bijection from the set
$S$ of  basis relations of $\cX$ onto the set of double 
$H$-cosets in~$G$.
Moreover, the intersection number $c_{rs}^t$ is equal to the 
multiplicity, with which an element of $D_t$ enters the product $D_r\,D_s$,
divided by $|H|$. It follows that
\qtnl{140714s}
n_s=\frac{|D_s|}{|H|}=\frac{|H|}{|H\cap H^g|}
\eqtn
for all $s\in S$ and $g\in D_s$ (the second equality follows
from the first one, because $|D_s|=|HgH|=|g^{-1}HgH|=|H^gH|$). 
In particular, 
$k=k(\cX)$ is the ratio between the order of $H$ and the minimal size 
of the intersection of $H$ with its conjugate.

\lmml{050216a}
Let $G$ be a transitive permutation group and $\cX=\inv(G)$.
If $b(\cX)\le 2$, then $G=\aut(\cX)$.
\elmm
\proof Inequality~\eqref{200615a} yields $b(G)\le b(\cX)\le 2$.
It follows that $H\cap H^g=1$ for some $g\in G$, where 
$H=G_\alpha$. If $s$ is the basis relation of $\cX$ with $D_s=HgH$, 
then~\eqref{140714s} implies 
that $\alpha s$ is a faithful regular orbit of~$H$. Hence
$$
|G|=nk,
$$ 
where $n$ is the cardinality of underlying set of $G$. 
Since $\cX=\inv(G)=\inv(\aut(\cX))$,
the above equality holds also for $G$ replaced by $\aut(\cX)$. Thus,
$$
|\aut(\cX)|=nk=|G|,
$$
and we are done, because $G\le\aut(\cX)$.\bull

\sbsnt{Indistinguishing number.}
Following \cite{MP12a}, the sum of all intersection numbers 
$c_{ss^*}^r$ with fixed $r$ is called the {\it indistinguishing number}
of $r\in S$ and denoted by $c(r)$. It is easily seen that for all
pairs $(\alpha,\beta)\in r$, we have
\qtnl{100814c}
c(r)=|\Omega_{\alpha,\beta}|,\ \text{where}\ 
\Omega_{\alpha,\beta}=\{\gamma\in\Omega:\ r(\gamma,\alpha)=r(\gamma,\beta)\}.
\eqtn
The maximum of the numbers $c(r)$, $r\ne 1_\Omega$, is called the
{\it indistinguishing number} of the coherent configuration~$\cX$
and denoted by $c(\cX)$.\medskip

The following lemma gives a formula for the indistinguishing number 
of the coherent configuration of a transitive permutation group. 
Recall that the fixity $\fix(G)$ of a permutation group $G$
is the maximum number of elements fixed by non-identity 
permutations \cite{SSh}.

\lmml{100814a}
Let $G\le\sym(\Omega)$ be a transitive group, $H$ a point stabilizer
of $G$, and $\cX=\inv(G)$. Then
\qtnl{100814d}
c(\cX)=\max_{x\in G\setminus H}|\bigcup_{h\in H}\Fix(hx)|.
\eqtn
In particular,
\qtnl{230415a}
c(\cX)\le  \max_{x\in G\setminus H}\sum_{h\in H}\chi(hx)\le \fix(G)\cdot|H|.
\eqtn
\elmm
\proof Let $r\in S$ and $(\alpha,\beta)\in r$. Then a point 
$\gamma$ belongs to the set $\Omega_{\alpha,\beta}$ defined 
in~\eqref{100814c} if and only if the pairs $(\gamma,\alpha)$
and $(\gamma,\beta)$ belong to the same orbit of the group~$G$ acting on $\Omega\times\Omega$, and the latter happens
if and only if $\gamma$ is a fixed point of a permutation $x\in G$ moving $\alpha$ to $\beta$.
Assuming without loss of generality that  $H=G_\alpha$, we conclude that  the set of all such $x$ forms an
$H$-coset $C$.  Therefore,
\qtnl{100814b}
c(r)=|\Omega_{\alpha,\beta}|=|\bigcup_{h\in H}\Fix(hx)|
\eqtn
for any $x\in C$. Moreover, if $r\ne 1_\Omega$, then $C\ne H$. This proves
equality \eqref{100814d}. Furthermore, $|\Fix(x)|=\chi(x)\le \fix(G)$ for any non-identity element $x\in G$.
This implies that
$$
|\bigcup_{h\in H}\Fix(hx)|\le
\sum_{h\in H}\chi(hx)\le \fix(G)\cdot|H|.
$$
Thus the second statement of the lemma follows from the first one.\bull

We complete this subsection by a statement that helps to 
count the values of the permutation character of a transitive group.

\lmml{100814e}
Let $G\le\sym(\Omega)$ be a transitive group, $\alpha\in\Omega$, 
and $H=G_\alpha$ the point stabilizer of $\alpha$ in~$G$. Then for 
every $x\in G$
\qtnl{080715a}
\Fix(x)\ne\varnothing\quad\Longleftrightarrow\quad x^G\cap H\neq\varnothing.
\eqtn
Suppose, additionally, that there is a subgroup $N$ with 
$H\le N\le N_G(H)$ such that every two elements of $H$ conjugated 
in $G$ are conjugated in $N$.
If $x=h_0^{g_0}$, where $h_0\in H$ and $g_0\in G$, then
\qtnl{260116a}
\Fix(x)=\{\alpha^g\mid g\in NCg_0\},
\eqtn
where $C=C_G(h_0)$. Furthermore,
\qtnl{080715b}
\chi(x)=\frac{|N:(C\cap N)|\,|C|}{|H|}=\frac{|N:(C\cap N)|\,|\Omega|}{|x^G|}.
\eqtn
\elmm

\proof  Clearly, $\alpha^g\in\Fix(x)$ if and only if $Hgx=Hg$, 
which holds if and only if there is $h\in H$ satisfying $x=h^{g}$. In particular, this yields~(\ref{080715a}).\medskip

To prove that the left-hand side of~\eqref{260116a} is contained
in the right-hand side, let $x=h_0^{g_0}$, that is the set $\Fix(x)$ is
nonempty. Suppose that $g$ is an arbitrary element of $G$ 
with $\alpha^g\in\Fix(x)$. Then there is $h\in H$ such that
$h^g=x=h_0^{g_0}$. Put $y=gg_0^{-1}$. Since the elements $h_0$ 
and $h=h_0^{y^{-1}}$ are conjugated in $G$, they are conjugated 
in~$N$, so there is  $n\in N$ with $h_0^{{y}^{-1}}=h_0^{{n}^{-1}}$. 
It follows that $y=nc$, where $c\in C$. Therefore, $g=ncg_0$, 
so $\alpha^g\in\Fix(x)$ implies that $g\in NCg_0$. To
establish the converse inclusion, for every $n\in N$, set
$h=h_0^{n^{-1}}$. Then $h^{ncg_0}=h_0^{cg_0}=x$ for
every $c\in C$. By the argument of the first paragraph, 
this proves $\alpha^{NCg_0}\subseteq\Fix(X)$.\medskip

Obviously, $|NCg_0|=|N:(C\cap N)||C|$. Now, the first equality
in~\eqref{080715b} is the direct consequence of~(\ref{260116a}),
because $\alpha^g=\alpha^{g'}$ if and only if $g'g^{-1}\in H$. 
Since $|C|=|G|/|x^G|$ and $|G|/|H|=|\Omega|$, the second equality
follows. \bull

\sbsnt{Algebraic isomorphisms and $m$-dimensional intersection numbers.}
Let $\cX=(\Omega,S)$ and $\cX'=(\Omega',S')$ be coherent configurations.
A bijection $\varphi:S\to S',\ r\mapsto r'$ is called an 
{\it algebraic isomorphism} from~$\cX$ to~$\cX'$ if
\qtnl{f041103p1}
c_{r^{}s^{}}^{t^{}}=c_{r's'}^{t'},\qquad r,s,t\in S.
\eqtn
In this case, we say that $\cX$ and $\cX'$ are {\it algebraically isomorphic}. Each isomorphism~$f$ from~$\cX$ to~$\cX'$ naturally 
induces an algebraic isomorphism between these coherent configurations.
The set of all isomorphisms inducing the algebraic isomorphism~$\varphi$
is denoted by $\iso(\cX,\cX',\varphi)$. In particular,
$$
\iso(\cX,\cX,\id_S)=\aut(\cX)
$$
where $\id_S$ is the identity on $S$. A coherent configuration~$\cX$ 
is called {\it separable} if, for any algebraic
isomorphism~$\varphi:\cX\to\cX'$, the set $\iso(\cX,\cX',\varphi)$ 
is nonempty.\medskip

Saying that coherent configurations $\cX$ and $\cX'$ have the same
intersection numbers, we mean that formula \eqref{f041103p1} holds
for a certain algebraic isomorphism. Thus, the exact meaning of 
the phrase "the coherent configuration $\cX$ is determined up 
to isomorphism by the intersection numbers" consists in the fact that
$\cX$ is separable.\medskip

Let $m\ge 1$ be an integer. According to~\cite{EP09},  the 
$m$-extension of a coherent configuration $\cX$ with point set 
$\Omega$ is defined to be the smallest coherent
configuration on $\Omega^m$, which contains the Cartesian $m$-power of
 $\cX$ and for which the set $\diag(\Omega^m)$ is a union of fibers. 
The intersection numbers of the $m$-extension are called the 
{\it $m$-dimensional numbers} of the configuration~$\cX$. Now,
$m$-separable coherent configurations for $m>1$ are defined 
essentially in the same way as for $m=1$. The exact definition
can be found in survey~\cite{EP09}, whereas in the present paper, we
need only the following result, which immediately follows 
from~\cite[Theorems~3.3 and~5.10]{EP09}.

\thrml{310116c}
Let $\cX$ be a coherent configuration admitting 
a $1$-regular  extension with respect to $m-1$ points, $m\ge 1$.
Then $\cX$ is $m$-separable.\bull 
\ethrm

\section{A sufficient condition for 1-regularity of a point
extension}\label{250415a}

\sbsnt{Main theorem.} The aim of this section is to prove the
following statement underlying the combinatorial part in the
proof of the main results of this paper.

\thrml{291114a}
Let $\cX$ be a homogeneous coherent configuration on $n$ points
with indistinguishing number $c$ and maximum valency $k$. Suppose
that $2c(k-1)<n$, i.e., inequality~\eqref{120415a} holds. Then every one point extension 
of $\cX$ is 1-regular. In particular, $b(\cX)\le 2$.
\ethrm

The proof of Theorem~\ref{291114a} will be given in the end of 
this section. The idea is to deduce the 1-regularity of the point extension $\cX_\alpha$
from Lemma~\ref{240714a} stating that inequality~\eqref{120415a}
implies the connectedness of the binary relations $s_{max}$ and
$s_\alpha$ defined in Subsection~\ref{020206a}. Note that this 
condition itself implies that any pair from $s_{max}$ forms a base 
of $\cX$ (Lemma~\ref{130714e}). 

\sbsnt{Relations $s_{max}$ and $s_\alpha$.}\label{020206a}
Let $\cX=(\Omega,S)$. Recall that $k=k(\cX)$ is the maximal valency of~$\cX$. Denote by $s_{max}$ the union of all relations in the set
$$
S_{max}=\{s\in S:\ n_s=k\}.
$$
Then, obviously, $S_{max}\subset S$ and $s_{max}\in S^\cup$. Moreover, since $\cX$ is homogeneous, we have $n_{s^*}=n_{s^{}}$ for all $s\in S$,
and hence, the relation $s_{max}$ is symmetric. We are interested in the connectedness of it, i.e., the connectedness
of the graph with vertex set $\Omega$ and edge set $s_{max}$. Note that, in general, this graph is not connected: take $\cX$ to be the
homogeneous coherent configuration of rank~$4$ that is associated
with a finite projective plane.\medskip

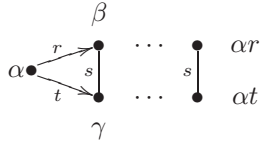
\begin{figure}[h]
	$\xymatrix@R=5pt@C=10pt@M=0pt@L=2pt{
		&   & \beta                &         &                       & \\
		&   & \bullet\ar@{-}[dd]_s & \cdots  & \bullet\ar@{-}[dd]_s  & \alpha r \\
		\alpha\bullet\ar@{->}[rrd]_t\ar@{->}[rru]^r &   &                      &         &                       & \\
		&   & \bullet              &  \cdots & \bullet               & \alpha t \\
		&   & \gamma               &         &                       & \\
	}$
	\caption{A part of the relation $s_\alpha$.}\label{lsp}
\end{figure}

With any point $\alpha\in\Omega$, we associate a binary relation $s_\alpha\subseteq \alpha s_{max}\times \alpha s_{max}$
that consists of all pairs $(\beta,\gamma)$ such that the colored triangle
$\{\alpha,\beta,\gamma\}$ is uniquely determined by the side colors $r=r(\alpha,\beta)$,
$s=r(\beta,\gamma)$  and $t=r(\alpha,\gamma)$, and one of the sides $\{\alpha,\beta\}$ or
$\{\alpha,\gamma\}$, see Fig.~\ref{lsp}. More precisely,
$$
s_\alpha=\{(\beta,\gamma)\in \alpha s_{max}\times \alpha s_{max}:\ c_{rs}^t=1\},
$$
This relation is symmetric. Indeed, we have $n_t=n_r=k$. Since also
$n_{r^*}=n_{r^{}}$, it follows from~\eqref{150410a} that $n_{t^{}}c_{r^{}s^{}}^{t^{}}=n_{r^*}c_{s^{}t^*}^{r^*}=n_{r^{}}c_{t^{}s^*}^{r^{}}$.
This implies that $c_{t^{}s^*}^{r^{}}=c_{r^{}s^{}}^{t^{}}=1$, and hence, $(\gamma,\beta)\in s_\alpha$.

\lmml{130216a}
Suppose that  the graph $s_\alpha$ is connected. Denote by 
$T_\alpha$ the set of all basis relations of the coherent
configuration $\cX_{\alpha}$ that are contained in 
$\alpha s_{max}\times \alpha s_{max}$. Then 
\qtnl{031015c}
|\beta t|=1\quad\text{for all}\ \,
t\in T_\alpha,\ \beta\in \alpha s_{max}.
\eqtn 
\elmm
\proof One can see that the set $\alpha s_{max}$ is the union of 
fibers of $\cX_\alpha$
(see~\cite[Lemma~2.2]{P12}). Therefore, 
\qtnl{210216a}
\alpha s_{max}\times \alpha s_{max}=\bigcup_{t\in T_\alpha}t.
\eqtn
Let  $t\in T_\alpha$ and $\beta\in\alpha s_{max}$. Then 
$\beta\in\alpha r$ for some $r\in S_{max}$. In view of~\eqref{210216a},  
there exists a point $\beta'\in\alpha s_{max}$ contained
in~$\beta t$.  By the connectedness
of $s_\alpha$, there exists a path~$P$ in $s_\alpha$ connecting 
$\beta$ and $\beta'$. If this path has length~$l=1$, then by
the definition of $s_\alpha$, we have
$c_{rt'}^s=1$, where $t'\in S$ and $s\in S_{max}$ are
unique relations such that $t\subset t'$ and $\beta'\in\alpha s$,
respectively. Then, obviously, $\{\beta'\}=\beta t$, 
as required.\medskip

Suppose that $l\ge 2$. Note that if $\beta_1$,
$\beta_2$, and $\beta_3$ are successive vertices of $P$, then they 
belong to $\alpha s_{max}$ and
$$
\{\beta_2\}=\beta_1t_1,\ \{\beta_3\}=\beta_2t_2,
$$
where $t_1$ and $t_2$ are the basis relations of $\cX_\alpha$ that
contain $(\beta_1,\beta_2)$ and $(\beta_2,\beta_3)$, respectively.
In particular, $t_1,t_2\in T_\alpha$ and 
$$
\{\beta_3\}=\beta_1 t_3,
$$
where $t_3$ is a unique relation of $T_\alpha$ containing the pair
$(\beta_1,\beta_3)$. This proves the required statement for $l=2$,
and, hence, for all positive integers~$l$ by induction.\bull

\lmml{130714e}
If $s_{max}$ and all $s_\alpha$, $\alpha\in\Omega$, are connected relations, then 
$\{\alpha,\beta\}$ is a base of the coherent configuration~$\cX$ for each
$\beta\in\Omega$ such that  $(\alpha,\beta)\in s_{max}$.
\elmm
\proof Let $\alpha\in\Omega$ and $\beta\in\alpha s_{max}$. Denote by
$\Gamma$ the set of all points $\gamma\in\Omega$ for which the singleton $\{\gamma\}$
is a fiber of the coherent configuration $\cX_{\alpha,\beta}$.
Then obviously $\alpha,\beta\in\Gamma$. We claim that
\qtnl{250511a}
\gamma s_{max}\subseteq\Gamma\qoq \gamma s_{max}\cap\Gamma=\varnothing
\eqtn
for all  $\gamma\in\Gamma$. Indeed, suppose on the contrary that there exist points $\gamma\in\Gamma$ and
$\gamma_1,\gamma_2\in\gamma s_{max}$ such that $\gamma_1\in\Gamma$ and $\gamma_2\not\in\Gamma$.
Since $s_\gamma$ is a connected relation, there is an $s_\gamma$-path connecting
$\gamma_1$ and $\gamma_2$. Moreover, the definition of $s_\gamma$ implies that if
some point in this path is inside $\Gamma$, then the next point in this path must
be also inside $\Gamma$. Therefore $\gamma_2\in\Gamma$,
a contradiction.\medskip

Denote by $\Gamma_0$ the set of all points $\gamma\in\Gamma$ with
$\gamma s_{max}\subseteq\Gamma$. Then $\alpha\in\Gamma_0$, because, as it follows from~\eqref{250511a}, the set $\alpha s_{max}$ contains $\beta\in\Gamma$.
Therefore, $\Gamma_0$ contains the connected component of  $s_{max}$ that
contains~$\alpha$. Since $s_{max}$ is connected, this implies that $\Gamma_0=\Omega$,
and hence $\Gamma=\Omega$. By the definition of $\Gamma$, this means that any fiber of the
coherent configuration $\cX_{\alpha,\beta}$ is a singleton. Thus,
$\{\alpha,\beta\}$ is a base of~$\cX$.\bull

\sbsnt{Connected components of $s_\alpha$.}
One can treat $s_\alpha$ also as the graph with vertex set 
$\alpha s_{max}$
and edge set $s_\alpha$. The set of all connected components of this graph that contain a vertex in $\alpha u$
for fixed $u\in S_{max}$ is denoted by $C_\alpha(u)=C(u)$.

\lmml{160814a}
Let $u,v\in S_{max}$. Suppose that $C(u)\cap C(v)\ne\varnothing$.
Then $C(u)=C(v)$ and $|\alpha u\cap C|=|\alpha v\cap C|$ for all
$C\in C(u)$.
\elmm
\proof Let $C_0\in C(u)\cap C(v)$. Then $C_0$ contains vertices
$\beta\in\alpha u$ and $\gamma\in\alpha v$
connected by an $s_\alpha$-path, say
$$
\beta=\beta_0,\beta_1,\ldots,\beta_m=\gamma,
$$
where $(\beta_i,\beta_{i+1})\in s_\alpha$ for~$i=0,\ldots,m-1$. By the definition of $s_\alpha$, this
implies that
\qtnl{291114w}
c_{u_iv_i}^{u_{i+1}}=1
\eqtn
for all $i$, where $u_i=r(\alpha,\beta_i)$ and  $v_i=r(\beta_i,\beta_{i+1})$. Therefore,
it is easily seen that for every $C\in C(u)$ given a vertex $\beta'\in C$, there is a unique $s_\alpha$-path
$$
\beta'=\beta'_0,\beta'_1,\ldots,\beta'_m=\gamma'
$$
such that $\gamma'\in \alpha v$ and $r(\alpha,\beta'_i)=u_i$ and  $r(\beta'_i,\beta'_{i+1})=v_i$ for all $i$. In
view of~\eqref{291114w}, no vertices $\beta^{}_i$ and $\beta'_i$ coincide whenever
$\beta\ne\beta'$. Thus, the mapping
$$
\alpha u\to\alpha v,\ \beta'\mapsto\gamma'
$$
is a bijection. Obviously, the vertex $\gamma'$ belongs to the component $C$
of the graph $s_\alpha$ that contains~$\beta'$. Since this is true for all $\beta'\in C$ and
all $C\in C(u)$, the required statement follows.\bull

For a relation $u\in S_{max}$ and a point $\delta\in\Omega$, denote by $p_u(\delta)$
the number of pairs $(\beta,\gamma)\in \alpha u\times\alpha u$ such that $\beta\ne\gamma$ and
$r(\beta,\delta)=r(\gamma,\delta)$. Here, $|\alpha u|=n_u=k$. Therefore, $\alpha u$ contains
exactly $k(k-1)$ pairs of distinct elements. Now we are able to estimate from above the sum of $p_u(\delta)$ in terms of the indistinguishing numbers of the corresponding basis relations $c(r(\beta,\gamma))$ as well as the indistinguishing number $c$ of $\cX$. Indeed,
\qtnl{301114z}
k(k-1)c\ge \sum_{\beta,\gamma}c(r(\beta,\gamma))\ge \sum_{\delta\in\Delta}p_u(\delta)
\eqtn
for any set $\Delta\subseteq\Omega$.
On the other hand, the number $p_u(\delta)$ can be computed by means of the intersection numbers.
Namely, if $v=r(\alpha,\delta)$, then, obviously,
\qtnl{160814b}
p_u(\delta)=\sum_{w\in T_{u,v}}c_{uw}^v(c_{uw}^v-1)
\eqtn
where  $T_{u,v}=\{w\in u^*v:\ c_{uw}^v>1\}$
(see Fig~\ref{f1}). In particular, the number $p_u(\delta)$ does not depend on $\delta\in\alpha v$.

\def\VRT#1{*=<5mm>[o][F-]{#1}}
\begin{figure}[h]
$\xymatrix@R=10pt@C=10pt@M=0pt@L=5pt{
&\VRT{\beta}\ar[drrrr]_{w} & & & &  \\
\VRT{\gamma}\ar[rrrrr]_{w} & & & & &\VRT{\delta}\\
 & & & & &\\
& & & & &\\
& & & \VRT{\alpha}\ar[llluuu]^u\ar[lluuuu]^u \ar[rruuu]_v& &  \\
}$
\caption{A relation $w\in T_{u,v}$ with $c_{uv}^w>1$.}\label{f1}
\end{figure}
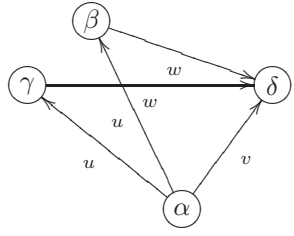

\lmml{160814c}
In the above notation, the following statements hold:
\nmrt
\tm{1} if either $n_u>n_v$, or $n_u=n_v$ and $C(u)\ne C(v)$, 
then $p_u(\delta)\ge k$,
\tm{2} if $n_u=n_v$, $C(u)=C(v)$, and $|C(u)|>1$, then $p_u(\delta)\ge k/2$.
\enmrt
\elmm
\proof To prove statement (1), suppose that either $n_u>n_v$, 
or $n_u=n_v$ and $C(u)\ne C(v)$. Then
\qtnl{160814d}
T_{u,v}=u^*v.
\eqtn
Indeed, obviously, $T_{u,v}\subseteq u^*v$. The converse inclusion is true if $n_u>n_v$, because in this case,
$c_{u^{}w^{}}^{v^{}}=n_uc_{v^*u^{}}^{w^*}/n_v>1$ for all $w\in u^*v$. Let now $C(u)\ne C(v)$.
Then the sets $C(u)$ and $C(v)$  are disjoint (Lemma~\ref{160814a}). This implies that
if $\beta\in \alpha u$ and $\gamma\in\alpha v$, then $(\beta,\gamma)\not\in s_\alpha$. Therefore,
$c_{uw}^v>1$ for all $w\in u^*v$, whence again $u^*v\subseteq T_{u,v}$. Thus relation~\eqref{160814d}
is completely proved. Together with \eqref{160814b}, this shows that
$$
p_u(\delta)=\sum_{w\in T_{u,v}}c_{uw}^v(c_{uw}^v-1)
\ge \sum_{w\in T_{u,v}}c_{uw}^v= \sum_{w\in u^*v}c_{uw}^v=n_u=k,
$$
as required. Observe that the penultimate equality is the well-known identity for homogenous coherent configurations. \medskip

To prove statement~(2), suppose that $n_u=n_v$, $C(u)=C(v)$, and $|C(u)|>1$. Let us choose $C\in C(u)$
so that the number $|\alpha u\cap C|$ is minimum possible. Then
\qtnl{301114a}
|\alpha u\setminus C|\ge k/2,
\eqtn
because $|C(u)|>1$ and $|\alpha u|=n_u=k$. Next, since $C(u)=C(v)$, we have $C\in C(v)$. Moreover,
$\alpha v$ is not contained in $C$, because $|C(v)|=|C(u)|>1$. Since $p_u(\delta)$ does not depend
on the choice of $\delta\in\alpha v$, we may assume that $\delta\in \alpha v\cap C$. Then  no point
$\beta\in \alpha u\setminus C$ belongs to the component of $s_\alpha$ that contains~$\delta$.
In particular, $(\delta,\beta)$ is not an edge of $s_\alpha$. Therefore,
$$
c_{uw}^v>1\quad\text{for all}\ w\in T,
$$
where $T$ is  the set of all  $w=r(\beta,\delta)$ with $\beta\in\alpha u\setminus C$.  By \eqref{160814b}
and \eqref{301114a}, we obtain
$$
p_u(\delta)=\sum_{w\in T_{u,v}}c_{uw}^v(c_{uw}^v-1)\ge
\sum_{w\in T}c_{uw}^v=|\alpha u\cap \delta T^*|\ge
|\alpha u\setminus C|\ge k/2,
$$
as required.\bull

\sbsnt{The connectedness of $s_{max}$ and $s_\alpha$.}
Using Lemmas \ref{160814a} and \ref{160814c}, we will prove that 
the hypothesis of Theorem~\ref{291114a} gives a sufficient condition
for the graphs $s_\alpha$ and $s_{max}$ to be connected. Note
that by Lemma~\ref{130714e}, this establishes the second statement
of Theorem~\ref{291114a}.

\lmml{240714a}
Suppose that $2c(k-1)<n$ and $k\ge 2$.  Then the graphs $s_\alpha$ and $s_{max}$ are connected.
Moreover, $|\alpha s_{max}|>n/2$.
\elmm
\proof To prove the first statement, we claim that
\qtnl{270116a}
|C(u)|=1\quad\text{for all}\ \, u\in S_{max}.
\eqtn
Indeed, if this is not true, then there exists $u\in S_{max}$ such that $|C(u)|\ge 2$. Lemma~\ref{160814c} yields that 
$p_u(\delta)\ge k/2$ for all points $\delta\in\Omega$. By~\eqref{301114z} with $\Delta=\Omega$, this implies that
$$
c\ge \frac{1}{k(k-1)}\sum_{\delta\in\Omega}p_u(\delta)\ge \frac{1}{k(k-1)}\frac{|\Omega|k}{2}
=\frac{n}{2(k-1)},
$$
which contradicts the lemma hypothesis. Thus, formula~\eqref{270116a} 
is proved.\medskip

Suppose on the contrary that the graph $s_\alpha$ is not connected for some $\alpha\in\Omega$. Then it has a component $C$
containing at most half of the vertices, that is
\qtnl{011214a}
2|C|\le |\alpha s_{max}|< n.
\eqtn
By \eqref{270116a}, one can find a relation $u\in S_{max}$ such that $C(u)=C$. Then for any point $\delta\in\Omega\setminus C$,
we have $n_v<n_u$ or $C(v)\ne C(u)$, where $v=r(\alpha,\delta)$ (if $C(v)=C(u)=C$, then $\delta\in C$).
By statement~(1) of Lemma~\ref{160814c}, this implies that $p_u(\delta)\ge k$. On the other hand,
$2|\Omega\setminus C|\ge n$ by \eqref{011214a}. From \eqref{301114z}  with
$\Delta=\Omega\setminus C$, we obtain that
\qtnl{011214e}
c\ge \frac{1}{k(k-1)}\sum_{\delta\in\Omega\setminus C}p_u(\delta)\ge
\frac{1}{k(k-1)}|\Omega\setminus C|k\ge \frac{n}{2(k-1)},
\eqtn
which contradicts the lemma hypothesis. Thus, the graph $s_\alpha$ is connected.\medskip

To prove that the graph $s_{max}$ is also connected, suppose on the contrary that one of its components,
say $C$, has at most $n/2$ points. Let $\alpha\in C$ and $u\in S_{max}$. Then
$\alpha u\subseteq C$ and $n_u>n_v$ for all $v=r(\alpha,\delta)$ with $\delta\in\Omega\setminus C$.
By statement~(1) of Lemma~\ref{160814c}, this implies that $p_u(\delta)\ge k$ 
for all such~$\delta$. Again inequality~\eqref{011214e} hold, which contradicts 
the lemma hypothesis.\medskip

To prove the second statement, denote by $V$ the union of all $v\in S$ with $n_v<k$, and fix $u\in S_{max}$. Then 
by statement (1) of Lemma~\ref{160814c}, we have $p_u(\delta)\ge k$ for all $\delta\in\Omega$ such that $r(\alpha,\delta)\in V$. 
This implies that
$$
\sum_{\delta\in \alpha V}p_u(\delta)\ge k|\alpha V|.
$$
Since there are $k(k-1)$ pairs of different points in $\alpha u$, at least for one of such pairs, say $(\alpha,\beta)$, we obtain
$$
c\ge c(u)=|\Omega_{\alpha,\beta}|\ge \frac{1}{k(k-1)}\sum_{\delta\in \alpha V}p_u(\delta)\ge \frac{k\,|\alpha V|}{k(k-1)}=
\frac{|\alpha V|}{k-1}.
$$
By the lemma hypothesis, this implies that $2\frac{|\alpha V|}{k-1}(k-1)<n$, whence $|\alpha V|<n/2$. Since $|\alpha V|+|\alpha s_{max}|=n$,
we are done.
\bull

\sbsnt{Proof of Theorem~\ref{291114a}.} Obviously, we may assume 
that $k\ge 2$. Let $\alpha\in\Omega$. It suffices to prove that 
each $\beta\in\alpha s_{max}$ is a regular  point of the coherent
configuration~$\cX_\alpha$. Suppose on the contrary that $\beta$ 
is not a regular point of
$\cX_\alpha$. Then this coherent configuration has a basis
relation~$v$ such that $|\beta v|\ge 2$. Therefore, there 
exist distinct points $\gamma_1$ and $\gamma_2$ such that
\qtnl{031015d}
r(\beta,\gamma_1)=v=r(\beta,\gamma_2).
\eqtn
Let $T_\alpha$ be the set defined in Lemma~\ref{130216a}. Then
in view of~\eqref{031015c}, $v\not\in T_\alpha$. Therefore,
neither $\gamma_1$ nor $\gamma_2$ belongs to $\alpha s_{max}$.\medskip

Let us verify that formula \eqref{031015d} holds with $\beta$
replaced by arbitrary $\beta'\in \alpha s_{max}$ and suitable 
distinct points $\gamma_1$ and $\gamma_2$. By virtue
of~\eqref{031015c}, the relation $u\in T_\alpha$ containing the 
pair $(\beta,\beta')$ is of the form
$$
u=\{(\delta,\delta^g):\ \delta\in\Delta \}
$$
for some bijection $g:\Delta\to\Delta'$, where $\Delta$ and 
$\Delta'$ are the fibers of $\cX_\alpha$ containing $\beta$ and
$\beta'$, respectively. Define a permutation
$f\in\sym(\Omega)$ by
$$
\omega^f=\css
\omega^g &\text{if $\omega\in  \Delta$,}\\
\omega^{g^{-1}} &\text{if $\omega\in  \Delta'$,}\\
\omega           &\text{otherwise}.\\
\ecss
$$
Then the graph of $f$ is the union of basis relations of $\cX_\alpha$, each of
which is of valency~$1$. One can see that in this case, $f$ takes any basis 
relation of $\cX_\alpha$  to another basis relation. This implies that $f$ is an isomorphism of the coherent 
configuration~$\cX_\alpha$ to itself. Therefore, 
by equality~\eqref{031015d}, we have
$$
r(\beta',\gamma_1)=r(\beta^f,\gamma_1^f)=
r(\beta^f,\gamma_2^f)=r(\beta',\gamma_2),
$$
and  the claim is proved. Thus, the set $\Omega_{\gamma_1,\gamma_2}$
contains $\alpha s_{max}$. By the theorem hypothesis and second statement of Lemma~\ref{240714a}, this implies that
$$
c\ge |\Omega_{\gamma_1,\gamma_2}|\ge |\alpha S_{max}|>n/2.
$$
However, then $n>2c(k-1)>n(k-1)$, which is impossible for $k>1$.\bull

\section{Inequality~\eqref{120415a} in simple groups of 
Lie type }\label{310116a}

The main purpose of the following two sections is to prove
Theorem~\ref{310116b} below, from which Theorems~\ref{030415a} and~\ref{300116c} were deduced in Introduction. In this section,
we reduce the proof to Lemma~\ref{030415b}, which will be proved
in the next section.

\thrml{310116b}
For the Cartan scheme $\cX$ of every group~$G\in\cL$,
inequality~\eqref{120415a} holds.
\ethrm

We generally follow
the notation of well-known Carter's book \cite{Car72} with some
exceptions that we explain below inside parentheses. If $\Phi_l$ 
is a simple Lie algebra of rank $l$, then $\Phi_l(q)$ is the simple
Chevalley group of rank $l$ over a field of order~$q$.  Let $B$, $N$,
and $H=B\cap N$ be a Borel, monomial, and Cartan subgroups of a
simple Chevalley group $\Phi_l(q)$ as in \cite{Car72}, while $W=N/H$
be the corresponding Weil group. Then \cite[Proposition 8.2.1]{Car72}
implies that the subgroups $B$ and $N$ form a $(B,N)$-pair 
of~$\Phi_l(q)$. If $\tau$ is a symmetry of the Dynkin diagram of
$\Phi_l$ of order $t$ (and the corresponding automorphism of
$\Phi_l(q^t)$), then ${}^t\Phi_l(q)$ is the simple twisted group of
Lie type (in \cite{Car72} such a group is denoted as ${}^t\Phi_l(q^t)$). Again $B$, $N$, and $H=B\cap N$ stand for Borel,
monomial, and Cartan subgroups of a simple twisted group of Lie type,
and $W=N/H$ is the Weil group (in \cite{Car72} they are denoted by
$B^1$, $N^1$ and so on). It follows from \cite[Theorem 13.5.4]{Car72}
that in this case $B$ and $N$ form a $(B,N)$-pair of~${}^t\Phi_l(q)$
again. In the sake of brevity we will use notation ${}^t\Phi_l(q)$
for all simple groups of Lie type, assuming that $t$ is the empty
symbol in the case of untwisted groups.  Recall also that the order
$w$ of the Weil group $W$ does not depend on the order of the
underlying field.\medskip

Let $G$ be a finite simple group of Lie type, and
let $\cX=(\Omega,S)$ be the Cartan scheme of $G$, where
the corresponding $(B,N)$-pair is as in the previous paragraph 
(see Definition~\ref{060216t}). In particular,
$\Omega=G/H$ and $S=\orb(G,\Omega^2)$. Put $n=|\Omega|$, 
$k=k(\cX)$, and $c=c(\cX)$.

\lmml{040216a}
There exists an element $g_0\in G$ such that $H\cap H^{g_0}=1$. In particular,
\qtnl{040216b}
k=\max_{s\in S}n_s=|H|.
\eqtn
\elmm

\proof It is well known that $B=U\rtimes H$ is the semidirect product
of $U$ and $H$, where $U$ is the unipotent radical of $B$. By
\cite[Propositions 5.1.5 and~5.1.7]{Car85}, there exists an unipotent
element $u\in U$ such that $C_G(u)\le U$ (in fact, $u$ can be chosen
as a regular unipotent element fixed by an appropriate Frobenius map
of the corresponding algebraic group), it follows that $[h,u]\neq1$
for every $h\in H^\#$. On the other hand, if $h\in H\cap H^u$, then
$[h,u]\in H\cap U=1$, so $h=1$ and $g_0=u$ is the desired element
of~$G$.\footnote{The alternative way to establish the same is to
apply Zenkov's theorem \cite{Ze}. It yields that since $H$ is
abelian, there is an element $g_0\in G$ such that $H\cap H^{g_0}$
lies in the Fitting subgroup of $G$, which is trivial in the case of
a simple group~$G$.} Now, in view of~\eqref{140714s}, 
$$
k=\max_{s\in S}n_s=\max_{g\in G}\frac{|H|}{|H\cap H^g|}=\frac{|H|}{|H\cap H^{g_0}|}=|H|.\,\bull
$$

Observe that $G$ satisfies the hypothesises of Lemmas~\ref{100814a}
and~\ref{100814e}. Indeed, the transitivity of the action of $G$ on
$\Omega$ is evident, while the monomial subgroup $N$ of $G$ satisfies
the additional condition from Lemma~\ref{100814e} due to
\cite[Proposition 8.4.5]{Car72}. This enables us to estimate the
indistinguishing number $c$ and prove the required
inequality~\eqref{120415a}. We need the following lemma whose
strictly group-theoretic proof is postponed to the next
section.\medskip

Below, for a coset  $\ov{y}=Hy$, set $M_{\ov{y}}=\{u\in \ov{y}:\  u^G\cap H\neq\varnothing\}$. Note that all elements of 
$M_{\ov{y}}$ are semisimple. For an integer $m$, set
$$
M_{\ov{y},m}=\{u\in M_{\ov{y}}:\ |u^G|\ge m\},\qaq M'_{\ov{y},m}=M_{\ov{y}}\setminus M_{\ov{y},m}.
$$
Put $$r_m=\max_{y\in G\setminus H}|M'_{\ov{y},m}|\qaq m_0=\min\limits_{\varnothing\neq y^G\cap H\neq y^G}|y^G|.$$

\lmml{030415b}
In the above notation, there exists a positive integer $m$ such that
\qtnl{030615a}
\frac{k}{m}+\frac{r_m}{m_0}\le\frac{1}{2wk}
\eqtn
for all groups $G\in\cL$ and every coset $\ov y\ne H$.
\elmm

{\bf Proof of Theorem~\ref{310116b}.} Immediately follows 
from Lemma~\ref{030415b} and Lemma~\ref{140216a} below.

\lmml{140216a}
Let $G$ be a simple group of Lie type. Suppose that there exists
an integer~$m$ such that inequality~\eqref{030615a} holds. Then 
for the Cartan scheme $\cX$ of~$G$, inequality~\eqref{120415a} is
satisfied.
\elmm
\proof It follows from Lemma~\ref{100814a} that
\qtnl{090715a}
c\le  \max_{y\in G\setminus H}\sum_{h\in H}\chi(hy)\le
\max_{y\in G\setminus H}\sum_{x\in M_{\ov{y}}}\chi(x).
\eqtn
Let $x\in M_{\ov y}$. Then by Lemma~\ref{100814e},
$\Fix(x)\ne\varnothing$, i.e. there are $h_0\in H$ and $g_0\in G$
with $x=h_0^{g_0}$. In the notation of that lemma 
$\chi(x)=|N:(C\cap N)|\,|\Omega|/|x^G|$. Furthermore, $|N:(C\cap N)|\le|N/H|=|W|=w$, because $H\le C\cap N$. We conclude that
\qtnl{080715c}
\chi(x)\le\frac{w\,n}{|x^G|}.
\eqtn

Let $m$ be a positive integer from Lemma~\ref{030415b}. Our
definitions imply that  $|x^G|\ge m$ for all $x\in M^{}_{\ov y,m}$
and  $|x^G|\ge m_0$ for all $x\in M_{\ov{y}}$. Taking into
account $|M_{\ov{y}}|\le|H|=k$ and formula~\eqref{080715c},
we obtain
$$
\sum_{x\in M_{\ov{y}}}\chi(x)\le
w\,n\sum_{x\in M_{\ov{y}}}\frac{1}{|x^G|}\le
w\,n\left(\sum_{x\in M^{}_{\ov y,m}}\frac{1}{|x^G|}+\sum_{x\in M'_{\ov y,m}}\frac{1}{|x^G|}\right)\le
$$
$$
w\,n\left(\frac{|M_{\ov y,m}|}{m}+\frac{|M'_{\ov y,m}|}{m_0}\right)\le
w\,n\left(\frac{k}{m}+\frac{r_m}{m_0}\right).
$$
By inequality~\eqref{030615a}, this implies
$$
\sum_{x\in M_{\ov{y}}}\chi(x)\le wn\cdot\frac{1}{2wk}=\frac{n}{2k}
$$
for any $y\in G\setminus H$. In view of~\eqref{090715a}, 
this immediately shows that $2ck\le n$, as required.\bull

\section{Proof of Lemma~\ref{030415b}}\label{070715a}

We begin with a simple remark. Suppose that $\ov{G}$ is a central
cover of $G$, i.e., $G\cong \ov{G}/Z(\ov{G})$, and subgroups
$\ov{B}$ and $\ov{N}$, whose images in $G$ are $B$ and $N$, form the
$(B,N)$-pair of $\ov{G}$. If $Z(\ov{G})\le\ov{H}=\ov{B}\cap\ov{N}$,
then we may exploit the action of $\ov{G}$ on its Cartan subgroup
$\ov{H}$ instead of the action of $G$ on $H$, because
$$
\inv(G,G/H)=\inv(\ov{G},\ov{G}/\ov{H}).
$$ 
Observe that the values of
$n$, $w$, the sizes of conjugacy classes remain the same, and $k=|H|$
does not exceed the order of $|\ov{H}|$. In particular, we may take
as $\ov{G}$ the universal cover $\widehat{G}$ of $G$ or, in the case
of classical groups, the group of linear transformation
$\widetilde{G}$ of vector space with an appropriate form, whose quotient by its central subgroup isomorphic to~$G$.\medskip

First, suppose that $G$ is a simple exceptional group. We prove
that relation~\eqref{030615a} holds for $m=m_0$. Note that the size of the conjugacy class of a semisimple element
of $G$ can be estimated from below by means of results in \cite{D83,D84} (for all exceptional groups other than the Ree and
Suzuki groups, this was done in \cite{Vd}). The corresponding low bounds for $m_0$ are listed in the second column
of Table~\ref{t1}. The values from the third and forth columns are well-known. Using this table, one can easily check
that $m_0\ge 2w|\widehat{H}|^2$, whence
\qtnl{100715b}
m_0\ge 2wk^2,
\eqtn
which is equivalent to inequality~\eqref{030615a} for $m=m_0$, because $r_{m_0}=0$ in this case.\medskip

\begin{table}[bt]
	\caption{}\label{t1}
	\begin{center}
		\begin{tabular}{|c|c|c|c|}
			\hline
			${}^t\Phi_l$        & $m_0$        & $|\widehat{H}|$            & $|W|$    \\
			\hline
			$E_8$      & $q^{112}$      & $(q-1)^8$        & $2^{14}\cdot3^5\cdot5^2\cdot7$ \\
			\hline
			$E_7$      & $(1/2)q^{64}$  & $(q-1)^7$        & $2^{10}\cdot3^4\cdot5\cdot7$ \\
			\hline
			$E_6$      & $(1/3)q^{30}$  & $(q-1)^6$        & $2^{7}\cdot3^4\cdot5$ \\
			\hline
			${}^2E_6$  & $(1/3)q^{30}$  & $(q-1)^4(q+1)^2$ & $2^{7}\cdot3^2$ \\
			\hline
			$F_4$      & $q^{16}$       & $(q-1)^4$        & $2^{7}\cdot3^2$ \\
			\hline
			$G_2$      & $q^{4}(q^3-1)$ & $(q-1)^2$        & $2^{2}\cdot3$ \\
			\hline
			${}^3D_4$  & $q^{16}$       &  $(q-1)(q^3-1)$  &  $2^2\cdot 3$\\
			\hline
			${}^2F_4$  &  $q^{6}(q-1)(q^3+1)$  &  $(q-1)^2$       & $2^4$ \\
			\hline
			${}^2G_2$  &  $q^2(q^2+q+1)$  &  $q-1$           & $2$ \\
			\hline
			${}^2B_2$  &  $q^2(q-1)$    &   $q-1$          & $2$  \\
			\hline
		\end{tabular}
	\end{center}
\end{table}

Now $G$ is a simple classical group. Our main source to 
estimate the size of a conjugacy class of~$G$ is
\cite{B}.\footnote{It is worth noting that in `an asymptotical sense'
 the required lower bounds can be taken from~\cite[Lemma~3.4]{LS}. 
We chose to use results from the later paper \cite{B} in order to
obtain the numerical values of $l_0$ and $a$.} Let $V$ be a natural
module over a field $\mathbb{F}_{q^v}$, where $v=2$ in the case of
unitary groups and $v=1$ otherwise, such that $G\le PSL(V)$, and let
$\widetilde{G}$ be preimage of $G$ in $SL(V)$. If $x\in G$ and 
$X\le G$, then $\widetilde{x}$ and $\widetilde{X}$ are preimages of
$x$ and $X$ in $\widetilde{G}$. We also agree to fix the base of $V$
in such a way that the preimage $\widetilde{H}$ of the Cartan
subgroup $H$ consists of diagonal matrices. Following
\cite[Definition 3.16]{B}, for an element $x\in G$ we denote by
$\nu(x)$ the codimension of largest eigenspace of $\widetilde{x}$ on
$\overline{V}=V\otimes K$, where $K$ is algebraic closure of
$\mathbb{F}_{q}$. For elements $x$ conjugated to elements of $H$,
which are of the prime interest for our purposes, $\nu(x)$ is just
equal to the difference between the dimension of $V$ and the maximum
eigenvalue multiplicity of  the diagonal matrix $\widetilde{h}$ with
$x=h^g$. We gather in Table~\ref{t2} the lower bounds $m_0$ on the sizes of conjugacy classes from \cite[Table 3.7-3.9]{B} as well as  the 
numbers $|W|$ and upper bounds for $|\widetilde{H}|$.\footnote{It is
worth mentioning that despite 
\cite[Tables 3.7-3.9]{B} contain the bounds on the sizes of 
conjugacy classes in the group $\operatorname{Inndiag}(G)$ rather
than $G$ itself, we get the correct bounds, because for $h\in H$ we
obviously have $|G:C_G(h)|=|\operatorname{Inndiag}(G):C_{\operatorname{Inndiag}(G)}(h)|$.} 
This table also contains the numbers $l_0$ and $a$ defining the
class~$\cL$. We also suppose that $q$ is odd in the case of the groups
$B_l(q)$ due to the well-known isomorphism $B_l(q)\cong C_l(q)$ for 
even $q$.\medskip

\begin{table}[bt]\label{010216a}
	\caption{}\label{t2}
	\begin{center}
		\begin{tabular}{|c|c|c|c|c|c|c|c|}
			\hline
			${}^t\Phi_l$   & conditions              & $m_0$              & $|\widetilde{H}|$                            & $|W|$                               & $l_0$  & $a$ \\
			\hline
			$A_l$          &                         & $q^{2l}/2$         & $(q-1)^l$                                    & $(l+1)!$                                & $7$    & $4$ \\
			\hline
			${}^2A_l$      & $l$ odd                 & $\frac{q^{4l-3}}{2(q+1)}$  & $(q-1)^{\lfloor\frac{l+1}{2}\rfloor}(q+1)^{\lfloor\frac{l}{2}\rfloor}$ & $2^{\lfloor\frac{l+1}{2}\rfloor}\lfloor\frac{l+1}{2}\rfloor!$ & $6$    & $4$ \\
			\hline
			${}^2A_l$      & $l$ even                & $\frac{q^{2l+1}}{2(q+1)}$  & $(q-1)^{\lfloor\frac{l+1}{2}\rfloor}(q+1)^{\lfloor\frac{l}{2}\rfloor}$     & $2^{\lfloor\frac{l+1}{2}\rfloor}\lfloor\frac{l+1}{2}\rfloor!$     & $6$    & $4$ \\
			\hline
			$B_l$          & $\frac{l(q-1)}{2}$ odd  & $\frac{q^{4l-1}}{4(q+1)}$  & $\frac{(q-1)^l}{2}$                          & $2^ll!$                             & $4$    & $4$ \\
			\hline
			$B_l$          & $\frac{l(q-1)}{2}$ even & $\frac{q^{2l+1}}{4(q+1)}$  & $\frac{(q-1)^l}{2}$                          & $2^ll!$                             & $4$    & $4$ \\
			\hline
			$C_l$          &                         & $\frac{q^{4l-4}}{2}$       & $(q-1)^l$                                    & $2^ll!$                             & $3$    & $4$ \\
			\hline
			$D_l$          &                         & $\frac{q^{4l-3}}{4(q+1)}$  & $(q-1)^l$                                    & $2^{l-1}l!$                         & $4$    & $2$ \\
			\hline
			${}^2D_l$      &                         & $\frac{q^{4l-3}}{4(q+1)}$  & $(q-1)^{l-1}(q+1)$                           & $2^{l-1}l!$                         & $4$    & $2$ \\
			\hline
		\end{tabular}
	\end{center}
\end{table}

Now we are ready to define the number $m$ for simple classical groups. Denote by $m_1$ the low bound for $|x^G|$ with $\nu(x)\ge 2$ that was found in \cite[Tables 3.7-3.9]{B}; the relevant values of $m_1$ are collected in the third column of Table~\ref{t3}. Set
\qtnl{160715a}
m=\css
m_0, &\text{if $\nu(h)\ge 2$ for all $h\in H^\#$,}\\
m_1, &\text{otherwise.}\\
\ecss
\eqtn
Note that for any coset $\ov y=Hy$  distinct from $H$,
$M^{}_{\ov y,m}=\{x\in M_{\ov y}\mid \nu(x)\ge 2\}$ and 
$M'_{\ov y,m}$ consist of all elements $x\in M_{\ov y}$ such that $\nu(x)=1$. Recall that $r_m=\max_{y\in G\setminus H}|M'_{\ov y,m}|$.

%

\lmml{040615a}
In the above notation the following statements hold.
\nmrt
\tm{1} If $G$ is one of the groups $A_l$, ${}^2A_l$ with $l$ even,
and $B_l$ with $l(q-1)/2$ even, then the number $r_m$ does not 
exceed the number in the fourth column of the corresponding row of 
Table~\emph{\ref{t3}}.
\tm{2} If $G$ is one of the other simple classical groups, then
$\nu(h)\ge2$ for every $h\in H^\#$.
\enmrt
\elmm

\begin{table}[bt]
	\caption{}\label{t3}
	\begin{center}
		\begin{tabular}{|c|c|c|c|}
			\hline
			${}^t\Phi_l$        & conditions              & $m_1$                      &  $r_m$     \\
			\hline
			$A_l$               &                         & $\frac{q^{4(l-1)}}{2}$     &  $\frac{l(l+1)(q-1)^2}{2}-1$  \\
			\hline
			${}^2A_l$           & $l$ even                & $\frac{q^{4l-3}}{2(q+1)}$  &  $\frac{(l+1)(q+1)^2}{2}+q$   \\
			\hline
			$B_l$               & $\frac{l(q-1)}{2}$ even & $\frac{q^{4l-1}}{4(q+1)}$  &  $\frac{l(q-3)}{2}+1$    \\
			\hline
		\end{tabular}
	\end{center}
\end{table}

\proof It is well known and easily verified that the diagonal
subgroup of a classical matrix group contains an element $h$ 
with $\nu(h)=1$ only if $G$ is one of the groups in statement~(1).
Therefore,  we need only to estimate~$r_m$ in these cases.\medskip

We claim that in any case
\qtnl{160715d}
r_m\le|\{h\in H^\#:\ \nu(h)\le 2\}|=:u.
\eqtn
Indeed, let $y$ and $z=hy$ be distinct elements of $M'_{\ov{y},m}$. Then $\nu(y)=\nu(z)=1$.
Therefore, each of the matrices  $\wt{y}$ and $\wt{z}$ has an eigenvalue of multiplicity $\dim(V)-1$. Since
these matrices are conjugated to diagonal matrices, this implies that the matrix $\wt{h}=\wt{y}(\wt{z})^{-1}$ has
an eigenvalue of multiplicity at least $\dim(V)-2$. Therefore, $\nu(h)\le2$. Since this is true for all
$z\in M'_{\ov{x},m}$, we are done.\medskip

Let $G=A_l$.  Then the required statement immediately follows
from~\eqref{160715d} by a direct calculation of the number $u$, 
which is equal the number of diagonal matrices in $\SL(l+1,q)$ 
with at most $2$ distinct diagonal enties.\medskip

Suppose that $G={}^2A_l$ and $l$ is even. To check the required upper bound on the number $u$, we observe that a base of $V$ can be chosen so that any matrix $\wt{h}\in\wt{H}$ is of the form
$$
\wt{h}=\operatorname{diag}(\lambda_1,\ldots,\lambda_r,\lambda_0,\lambda_1^{-q},\ldots,\lambda_r^{-q}),
$$
where $r=l/2$, $\lambda_i\in \mathbb{F}_{q^2}$ for all $i$, $(\lambda_0)^{q+1}=1$, and
$\lambda_0(\lambda_1)^{1-q}\cdots(\lambda_r)^{1-q}=1$. If, in addition, $\nu(h)\le2$ and $l\ge6$, then either
$$
\widetilde{h}=\operatorname{diag}(\lambda,\ldots,\lambda,\lambda_0,\lambda,\ldots,\lambda),
$$
where $\lambda^{q+1}=1$ and $\lambda_0\lambda^{l}=1$, or
$$\widetilde{h}=\operatorname{diag}(\lambda,\ldots,\lambda,\mu,\lambda,\ldots,\lambda,\lambda_0,\lambda,\ldots,\lambda,\mu^{-q},\lambda,\ldots,\lambda),
$$
where $\lambda=\lambda_0$, $\lambda^{l-1}\mu^{1-q}=1$, and $\mu$ takes an arbitrary $j$-th of the first $r$ positions (so $\mu^{-q}$
takes the $(r+1+j)$-th position). The rest is routine.\medskip

Let $G=B_l$ and $l(q-1)/2$ even. Then $G=\wt{G}=\Omega_{2l+1}(q)$. 
To estimate $u$ from above, choose a base of $V$ so that any
matrix $h\in H$ is of the form
$$
h=\operatorname{diag}(\xi^{k_1},\ldots,\xi^{k_l},\xi^{-k_1},\ldots,\xi^{-k_l},1),
$$
where $\xi$ is a primitive element of the field $\mF_q$, and the number $k_1+\ldots+k_l$ is even. If, in addition, $\nu(h)\le2$
and $l\ge3$, then either
$$
h=\operatorname{diag}(-1,\ldots,-1,1)
$$
(recall that $l(q-1)/2$ is even), or
$$h=\operatorname{diag}(1,\ldots,1,\mu,1,\ldots,1,\mu^{-1},1,\ldots,1),
$$
where $\mu$ is a nonzero square in $\mF_q$ and takes an arbitrary $j$-th of the first $l$ positions (so $\mu^{-1}$ takes the $(l+j)$-th position). Thus, $u\le l(q-3)/2+1$.\bull\medskip

To complete the proof, we verify inequality~\eqref{030615a} for 
the number $m$ defined by~\eqref{160715a}. Observe that, due
to~\eqref{160715a} and Lemma~\ref{040615a}, the number $r_m$
equals $0$ in all cases when $m=m_0$. In the latter case, 
it suffices to verify inequality~\eqref{100715b}. 
We proceed further case by case.\medskip

Let $G=C_l(q)$. Here $m=m_0$ and we need to prove that 
$m_0\ge 2wk^2$.  According to Table~\ref{t2}, it means that for
$l\ge3$ and $q\ge4l$, the following inequality have to be true:
$$
\frac{q^{4l-4}}{2}\ge 2^{l+1}l!(q-1)^{2l}.
$$
For $l=3$, this inequality is straightforward. If $l\ge 4$, then $q^{2l-4}\ge 4(2l)^l\ge 2^{l+2}l!$, and we are done.\medskip

Let $G=D^{\varepsilon}_{l}(q)$. Then $m=m_0$ and to
verify~\eqref{100715b}, we check that
$$
\frac{q^{4l-3}}{4(q+1)}\ge 2^ll!(q-1)^{2l-2}(q+1)^2
$$
for $l\ge4$ and $q\ge2l$. Since $(q-1)^{2l-2}(q+1)^3<q^{2l+1}$ for all these $l$ and $q$, it suffices to check that $q^{2l-4}\ge 2^{l+2}l!$. For $l=4$ it can be verified directly, while for $l>4$ we have $q^{2l-4}\ge4(2l)^l\ge2^{l+2}l!$.\medskip

Let $G=A_{l}(q)$ and suppose that $l\ge 7$ and $q\ge4l$. By Lemma~\ref{040615a}, we obtain that
$$
r_m\le \frac{l(l+1)(q-1)^2}{2}-1\le \frac{q^4}{32}.
$$
Thus, the left-hand side of~\eqref{030615a} can be estimated as follows:
\qtnl{160715x}
\frac{k}{m}+\frac{r_m}{m_0}\le \frac{2(q-1)^l}{q^{4l-1}}+\frac{2q^4}{32q^{2l}}\le \frac{2}{q^{4l-3}}+\frac{1}{16q^{2l-4}}\le
\frac{1}{8q^{2l-4}}.
\eqtn
On the other hand, $(l+1)!\le 4^{l-3}l^{l-4}$ for $l\ge 7$: this 
is verified directly for $7\le l\le 9$, and follows from
the obvious inequalities $(l+1)!<l^l<4^{l-3}l^{l-4}$ for $l\ge 10$. Therefore, in our case, we get the following
low bound for the right-hand side of~\eqref{030615a}:
\qtnl{160715y}
\frac{1}{2wk}=\frac{1}{2(q-1)^l(l+1)!}\ge \frac{1}{2q^l 4^{l-3}l^{l-4}}\ge \frac{1}{2q^l 4q^{l-4}}=\frac{1}{8q^{2l-4}}.
\eqtn
Thus, the required statement follows from~\eqref{160715x} and~\eqref{160715y}.\medskip

For each of the remaining two series of classical groups the
expression on the left-hand side of~\eqref{030615a} for 
$m=m_0\le m_1$ does not exceed the
same expression for $m=m_1$ (see Tables~\ref{t2} and~\ref{t3}). Since the expression on the right-hand side in both cases does not 
depend on whether $m=m_0$ or not,
it suffices to verify~\eqref{030615a} for $G={}^2A_{l}(q)$ (resp.,
$G=B_{l}(q)$) independently of the oddness of $l$ (resp.,
$l(q-1)/2$), taking $m_0$ and $m$ as in case of even $l$ (resp. 
even $l(q-1)/2$).\medskip

Let $G={}^2A_{l}(q)$. Suppose that $l\ge 6$ and $q\ge4l$. Lemma~\ref{040615a} yields that
$$
r_m\le \frac{(l+1)(q+1)^2}{2}+q\le \frac{q^3}{6}.
$$
Put $b=\lfloor l+1/2\rfloor$. Now, the left-hand side and right-hand side of~\eqref{030615a} can be estimated as follows:
\qtnl{160715r}
\frac{k}{m}+\frac{r_m}{m_0}\le
\frac{2(q-1)^{b}(q+1)^{\lfloor\frac{l}{2}\rfloor+1}}{q^{4l-3}}+\frac{2q^3(q+1)}{6q^{2l+1}}\le
\frac{2q^l(q+1)}{q^{4l-3}}+\frac{q^3(q+1)}{3q^{2l+1}}
\eqtn
and
\qtnl{160715e}
\frac{1}{2wk}\ge
\frac{1}{2(q-1)^{b}(q+1)^{\frac{l}{2}}2^{b}b!}\ge
\frac{1}{2q^l2^bb!}.
\eqtn
By~\eqref{160715r} and~\eqref{160715e}, it suffices to verify that
\qtnl{160715q}
2^bb!\le\frac{3q^{2l-3}}{2(q+1)(q^{l-1}+6)}.
\eqtn
However, one can easily check that $2(q+1)(q^{l-1}+6)\le 3q^{l}$ and $2^bb!\le q^{l-3}$ for all
$q\ge 4l\ge 25$. Therefore, \eqref{160715q} holds, and we are done.\medskip

Let $G=B_{l}(q)$. Suppose that $l\ge 4$ and $q\ge 4l$. By Lemma~\ref{040615a}, it follows that
$$
r_m\le \frac{l(q-3)}{2}+1\le \frac{q^2}{8}.
$$
Now, the left-hand side and right-hand side of~\eqref{030615a} can be estimated as follows:
$$
\frac{k}{m}+\frac{r_m}{m_0}\le
\frac{4q^l(q+1)}{2q^{4l-1}}+\frac{4q^2(q+1)}{8q^{2l+1}}=
\frac{(q+1)(q^l+4)}{2q^{3l-1}}
$$
and
$$
\frac{1}{2wk}\ge \frac{1}{2^ll!(q-1)^l}\ge \frac{1}{2^ll!q^{l-1}(q-1)}.
$$
Thus, it suffices to verify that
$$
(q^l+4)2^ll!\le 2q^{2l-2}.
$$
This is straightforward for $l=4$. Since $q\ge 4l$, the required inequality holds whenever $l!\le 2^{l-4}l^{l-2}$, which can be
directly checked for $5\le l\le 10$. Finally, if $l\ge11$, then
 $l!\le l^l\le 2^{l-4}l^{l-2}$. This completes the proof of the
lemma.

\section{Proof of Theorem~\ref{120714a}}\label{240714z}
Note that the actions of $\PSL(2,q)$ and $G=\SL(2,q)$ on the cosets 
of the corresponding Cartan subgroups are equivalent. Thus, without
loss of generality, we may assume that $\cX=\inv(G,\Omega)$, where
$\Omega=G/H$ with $H$ being the subgroup of diagonal matrices 
of $G$. Thus, 
\qtnl{030216a}
|G|=q(q+1)(q-1),\qquad |H|=q-1,\qquad |\Omega|=|G:H|=q^2+q.
\eqtn
First, we study a structure of $\cX$ in terms of double
$H$-cosets (see Subsection~\ref{030216f}).\medskip

One can see that the group $N=N_G(H)$ is the disjoint union of
two double $H$-cosets, namely, $H$ and $HiH=Hi$, where 
$$
i=\begin{pmatrix}
0 & -1 \\
1 & \phantom{-}0 \\
\end{pmatrix}
$$
Denote by $s_1$ and $s_i$ the basis relations of $\cX$, for which
$D_{s_1}=H$ and $D_{s_i}=HiH$ (see~\eqref{140714u}). Clearly,
$s_1=1_\Omega$.

\lmml{030216b}
Let $S$ be the set of basis relations of the coherent
configuration~$\cX$. Then given $s\in S$, we have
$$
n_s=\css
1    &\text{if $s\in\{s_1,s_i\}$},\\
q-1  &\text{otherwise.}\\
\ecss
$$
In particular, $|S|=q+4$ and $|S_{max}|=q+2$.
\elmm
\proof It is easy to verify that $H^x\cap H=1$ for all
$x\in G\setminus N$ and $N=H\cup Hi$. Thus, the required statements
follow from formula~\eqref{140714s}.\bull


Denote by $U$ and $V$ the subgroups (in $G$) of unipotent upper
triangular and low triangular matrices, respectively. Since, obviously,
$H\le N_G(U)\cap N_G(V)$, we conclude that
\qtnl{030216h}
HuH=HU^\#=U^\#H\qaq HvH=HV^\#=V^\#H
\eqtn
for all $u\in U^\#$ and $v\in V^\#$. Denote by $s_u$ and $s_v$
the basis relations of $\cX$, for which $D_{s_u}=HuH$ and
$D_{s_v}=HvH$, respectively. In view of~\eqref{030216h},
these relations do not depend on the choice the matrices
$$
u=\begin{pmatrix}
1 & x \\
0 & 1 \\
\end{pmatrix},\qaq
v=\begin{pmatrix}
1 & 0 \\
y & 1 \\
\end{pmatrix},\quad
$$
where $x$ and $y$ are nonzero elements of the field~$\mathbb{F}_q$.
Clearly, from~\eqref{030216h}, it follows that $(s_u)^*=s_u$ and
$(s_v)^*=s_v$.

\lmml{190714e}
In the above notation, let $s\in S$. Then
\nmrt
\tm{1} $c_{s_u s_{}}^{s_v}=0$ if $s=s_1$ or $s_i$, and
$c_{s_u s_{}}^{s_v}=1$ otherwise,
\tm{2} if $s\not\in\{s_1,s_i,s_u,s_v\}$, then $c_{s_u^{}s^{}_v}^s=1$ or
$c_{s_v^{}s^{}_u}^s=1$.
\enmrt
\elmm
\proof It is straightforward to check that 
$c_{s_u s_1}^{s_v}=c_{s_u s_i}^{s_v}=0$. Due to Lemma~\ref{030216b},
we may assume that $s\in S_{max}$. Then by this lemma, 
$n_{s_u}=q-1=n_{s_v}$. Therefore, 
$c_{s_u s_{}}^{s_v}=c_{s_u s_{v}}^s$.
The number $|H|\,c_{s_u s_{v}}^s$ is equal to the multiplicity, with which
an element $w\in D_s$ enters the product
$$
D_{s_u}\,D_{s_v}=HuH\, HvH=HU^\#\,HV^\#=HH(U^\#V^\#)
$$ 
(see \eqref{030216h}). Thus, to prove statement~(1), it suffices to verify
that no two elements in $UV$ belong to the same $H$-coset. For this aim, 
suppose that  $u_1v_1h=u_2v_2$ for some $u_1,u_2\in U$, $h\in H$, and 
$v_1,v_2\in V$. Then the group $U$ of unipotent upper triangular
matrices contains the element
$$
u_2^{-1}u_1^{\phantom{-1}}=v_2^{}\,h_{}^{-1}v_1^{-1},
$$
which is a low triangular matrix. It follows that $u_1=u_2$, 
$v_1=v_2$, $h=1$, and  we are done.\medskip

To prove statement (2), it suffices to verify that the complement to
the set $D_{s_u}D_{s_v}\,\cup\, D_{s_v}D_{s_u}$ in $G$ is equal to
$D_{s_1}\,\cup\,D_{s_i}\,\cup\,D_{s_u}\,\cup\,D_{s_v}$. In view of
equalities~\eqref{030216h}, this is equivalent to
\qtnl{040216r}
G\setminus (HU^\#V^\#\,\cup\, HV^\#U^\#)=N\,\cup\,HU^\#\,\cup\, HV^\#.
\eqtn
To prove this relation, we observe that general elements of the sets 
$U^\#V^\#$ and $V^\#U^\#$ are, respectively, 
$$
\begin{pmatrix}
1+xy & x \\
y & 1 \\
\end{pmatrix}
\qaq
\begin{pmatrix}
1 & x \\
y & 1+xy \\
\end{pmatrix},
$$ 
where $x$ and $y$ are nonzero elements in $\mathbb{F}_q$. Therefore, there are 
at least $q-1$ elements in  $V^\#U^\#$, which do not belong to $U^\#V^\#$ 
(they correspond to  nonzero elements $xy$). By the statement proved in
the previous paragraph, it follows that the set
$HU^\#V^\#\,\cup\, HV^\#U^\#$ is the disjoint union of $(q-1)^2$ distinct
cosets of $H$ contained in  $HU^\#V^\#$ and at least $q-1$ distinct
cosets of $H$ contained in  $HV^\#U^\#$. Since none of all these
cosets is contained in $N\,\cup\,HU^\#\,\cup\, HV^\#$, we have
$$
(q-1)^2q=q(q+1)(q-1)-2(q-1)-2(q-1)^2=
$$
$$
|G\setminus(N\,\cup\,HU^\#\,\cup\, HV^\#)|\ge
|HU^\#V^\#\,\cup\, HV^\#U^\#|\ge 
$$
$$
(q-1)^2(q-1)+(q-1)(q-1)=(q-1)^2q,
$$
which proves formula~\eqref{040216r}.\bull

Let us verify that the coherent
configuration $\cX_\alpha$ with $\alpha=H$, is $1$-regular. Indeed,
in this case the first statement of Theorem~\ref{120714a}
follows from Theorem~\ref{310116c} for $m=2$, whereas the second
statement is obvious.\medskip 

To prove the $1$-regularity of $\cX_\alpha$, it suffices to check
that every point $\beta\in\alpha s_{max}$ is regular. However,
if $t$ is a basis relation of $\cX_\alpha$, then $t$ is contained
in a basis relation $s$ of~$\cX$. If $s\in\{s_1,s_i\}$, then by
Lemma~\ref{030216b}, we have
$$
|\beta t|\le|\beta s|=n_s=1.
$$
Thus, by the same lemma, we may assume that $s\in S_{max}$
and hence $t$ belongs to the set $T_\alpha$ defined in
Lemma~\ref{130216a}. By this lemma, all we need is to verify
the connectedness of the graph $s_\alpha$.\medskip

Let us prove the connectedness of $s_\alpha$. 
Suppose that the pair $(\gamma, \delta)\in\alpha u\times\alpha v$,
belongs to the basis relation~$s$. Then, obviously, 
$c_{s_u s_{}}^{s_v}\ne 0$. Thus, by statement~(1) of
Lemma~\ref{190714e} and the definition of $s_\alpha$, the points
$\gamma$ and $\delta$ are adjacent in $s_\alpha$.
From statement~(2) of Lemma~\ref{190714e}, it follows
that any other vertex 
$\beta\in\alpha s$ with $s\in S_{max}$, has at least one neighbor 
in the set 
$\alpha s_u\,\cup\,\alpha s_v$, i.e.,
$$
\beta s_\alpha\,\cap\, (\alpha s_u\,\cup\,\alpha s_v)\, \ne\, \varnothing.
$$
Thus, $s_\alpha$ is connected, and we are done.

\section{Proof of Theorem~\ref{310116e}}\label{210615b}

We make use of the well-known Weisfeiler-Leman algorithm described in detail
in~\cite[Section~B]{W76}. The input of it is a set $\cS$ of binary 
relations on a set $\Omega$, and the output is the smallest coherent
configuration 
$$
\WL(\cS)=(\Omega,S)
$$ 
such that  $\cS\subset S^\cup$. The running time of 
the algorithm is polynomial in the cardinalities of~$\cS$
and~$\Omega$. The proof of the following statement is based on
the Weisfeiler-Leman algorithm and can be found in \cite[Theorem~3.5]{P12}.

\thrml{270411s}
Let $\cX$ and $\cX'$ be coherent configurations on $n$ points. Then given 
an algebraic isomorphism $\varphi:\cX\to\cX'$ all the elements of the set $\iso(\cX,\cX',\varphi)$ can be listed in time $(bn)^{O(b)}$ where $b=b(\cX)$.
\ethrm

To solve the recognition problem, at first, we recognize the colored 
graphs $D$ of Cartan schemes of $G$ with respect to $(B,N)$-pair of 
rank at least~$2$. In this case, from the corollary of the main
theorem in~\cite{MHTM}, it follows that  $B$ is the normalizer
$N_G(P)$ of a group $P\in\syl_p(G)$ such that 
\qtnl{060216a}
H\cap P=1,\qquad H\le N_G(P),\qquad |N_G(P)|=|H|\,|P|,
\eqtn
where $p$ is the characteristic of the ground field. 
By \cite{V74}, one can also see that apart for a finite number of
exceptional groups, $N=N_G(H)$  for every group from~$\cL$.
Thus, the correctness of 
the following algorithm follows from Theorems~\ref{300116c}
and~\ref{270411s}.\medskip

In what follows, we denote by $\Omega$ the 
vertex set of the graph $D\in\cG_n$, by $\cS$ the set of its color
classes, and by $\cS_{\alpha,\beta}$ the 
union of $\cS$  and the set of two one-element relations
$\{(\alpha,\alpha)\}$ and $\{(\beta,\beta)\}$.\medskip

\centerline{\bf Recognizing Cartan schemes 
(the rank of $(B,N)$ is at least~2)}\medskip

\noindent{\bf Step 1.} Find the coherent configuration
$\cX=\WL(\cS)$.\vspace{1mm}

\noindent{\bf Step 2.} If there are no distinct points
$\alpha,\beta$ such that the coherent configuration
$\cX_{\alpha,\beta}=\WL(\cS_{\alpha,\beta})$ is complete,
then $b(\cX)>2$ and $D\not\in\cK_n$.\vspace{1mm}

\noindent{\bf Step 3.} Find all the elements of the group
$G=\iso(\cX,\cX,\id)$ by the algorithm of Theorem~\ref{270411s}.
If $G$ is not simple, then $D\not\in\cK_n$.\vspace{1mm}

\noindent{\bf Step 4.} Analyzing the number $|G|$, check that 
$G\in\cL$. If not, then $D\not\in\cK_n$; otherwise set $p$ to 
be the characteristic of the ground field associated with~$G$.
\vspace{1mm}

\noindent{\bf Step 5.}
Fix a point stabilizer $H$ of $G$ and find $P\in\syl_p(G)$,
for which relations~\eqref{060216a} hold. If there is no such~$P$,
then $D\not\in\cK_n$. \vspace{1mm}

\noindent{\bf Step 6.} Now $D\in\cK_n$ and $\cX$ is the Cartan
scheme of $G$ with respect to $(B,N)$,  where $B=N_G(P)$ and $N=N_G(H)$.\bull\medskip

Let us estimate the running time of the algorithm. At  Steps~1 
and~2, we apply the Weisfeler-Leman algorithm $n(n-1)+1$ times.
Thus, the complexity of these steps is at most $n^{O(1)}$. At
Step~3, the time is polynomially bounded by Theorem~~\ref{270411s}
and the fact that a group is simple if and  only if 
no nontrivial conjugacy class of it generates a proper subgroup
(given the elements of $G$ the conjugacy classes of it can be
found efficiently). Step~4 requires polynomially many of
arithmetic operations involving the number~$|G|$ written in unary
system. Here, we use the fact based on CFSG that except for 
known cases, any finite simple group is uniquely determined by 
its order (see Theorem~5.1 and Lemma~2.5 in~\cite{KLST}). 
Since Steps~5 and~6 can obviously be implemented in
polynomial time for the group~$G$ given by the multiplication table,
we conclude that the running time of the algorithm is 
at most~$n^{O(1)}$.\medskip

The first four steps of the algorithm remain the same as before
if if we do not assume that the rank of $(B,N)$ is at least~2.
But in this case, one can find a $2$-transitive representation
of the group $G$; here, a complete classification of all
$2$-transitive groups is useful (see, e.g., \cite{DM}).  This 
enables to find the group~$B$ and~$N$.\medskip

To solve the isomorphism problem, let $D\in\cK_n$ and $D'\in\cG_n$.
Denote by $\cS$ and $\cS'$ the sets of color classes of $D$ and $D'$,
respectively. Without loss of generality, we may assume that there
is a color preserving bijection $\psi:\cS\to\cS'$. Then one can
apply the canonical version of the Weisfeiler-Leman 
algorithm presented in~\cite[Section~M]{W76}, where, in fact, the following statement was proved.

\thrml{thbw}
Let $\cS$ and $\cS'$ be $m$-sets of binary relations on an $n$-element set.
Then given a bijection $\psi:\cS\to\cS'$ one can check in time $mn^{O(1)}$
whether or not there exists an algebraic isomorphism
$\varphi:\WL(\cS)\to\WL(\cS')$ such that $\varphi|_\cS=\psi$. Moreover, 
if $\varphi$ does exist, then it can be found within the same time.\bull
\ethrm

Clearly, the original graphs $D$ and $D'$ are not isomorphic if 
there is no algebraic isomorphism $\varphi$ from Theorem~\ref{thbw}.
Assuming the existence of $\varphi$, we can find the set 
$$
\iso(D,D')=\iso(\cX,\cX',\varphi)
$$
in time $(bn)^{O(b)}$ by Theorem~\ref{270411s}, where $\cX=\WL(\cS)$,
$\cX'=\WL(\cS')$, and $b=b(\cX)$. Since $b\le 2$ (Theorem~\ref{300116c}),
we are done.

\end{document}